\documentclass{article}

\newtheorem{theorem}{Theorem}[section]

\newtheorem{proposition}[theorem]{Proposition}
\newtheorem{definition}[theorem]{Definition}

 \evensidemargin\oddsidemargin
\begin{document}

\title{Quasi-concave functions on antimatroids}
\author{Vadim E. Levit and Yulia Kempner \\
%EndAName
Department of Computer Science\\
Holon Academic Institute of Technology\\
52 Golomb Str., P.O. Box 305\\
Holon 58102, ISRAEL \\
\{yuliak, levitv\}@hait.ac.il}
\maketitle

\begin{abstract}
In this paper we consider quasi-concave set functions defined on
antimatroids. There are many equivalent axiomatizations of antimatroids,
that may be separated into two categories: antimatroids defined as set
systems and antimatroids defined as languages. An algorthmic
characterization of antimatroids, that considers them as set systems, was
given in \cite{KL}. This characterization is based on the idea of
optimization using set functions defined as minimum values of linkages
between a set and the elements from the set complement. Such set functions
are quasi-concave. Their behavior on antimatroids was studied in \cite{Malta}%
, where they were applied to constraint clustering. In this work we
investigate a duality between quasi-concave set functions and linkage
functions. Our main finding is that an arbitrary quasi-concave set function
on antimatroid may be represented as minimum values of some monotone linkage
function.\newline

\textbf{keywords:}\textit{\ antimatroid, quasi-concave set function,
monotone linkage function.}
\end{abstract}

\section{ Introduction}

Let $E$ be a finite set. A \textit{set system }over $E$\textit{\ }is a pair $%
(E,\mathcal{F})$, where $\mathcal{F}\subseteq 2^{E}$ is a family of subsets
of $E$, called \textit{feasible }sets. We will use $X\cup x$ for $X\cup
\{x\} $, and $X-x$ for $X-\{x\}$.

\begin{definition}
A non-empty set system $(E,\mathcal{F})$ is an \textit{antimatroid }if

$(A1)$ for each non-empty $X\in \mathcal{F}$, there is an $x\in X$ such that 
$X-x\in \mathcal{F}$

$(A2)$ for all $X,Y\in \mathcal{F}$, and $X\not\subseteq Y$, there exist an $%
x\in X-Y$ such that $Y\cup x\in \mathcal{F}$.
\end{definition}

\smallskip Any set system satisfying $(A1)$ is called \textit{accessible.}

\begin{definition}
\smallskip A set system $(E,\mathcal{F})$ has the \textit{interval property
without upper bounds }if for all $X,Y\in \mathcal{F}$ with $X\subseteq Y$
and for all $x\in E-Y$ , $X\cup x\in \mathcal{F}$ implies $Y\cup x\in 
\mathcal{F}$.
\end{definition}

\smallskip There are some different antimatroid definitions:

\begin{proposition}
\cite{BZ}\cite{Greedoids}For an accessible set system $(E,\mathcal{F})$ the
following statements are equivalent:

$(i)$ $(E,\mathcal{F})$ is an antimatroid

$(ii)$ $\mathcal{F}$ is closed under union

$(iii)$ $(E,\mathcal{F})$ satisfies the interval property without upper
bounds.
\end{proposition}

For a set $X\in \mathcal{F}$, let $\Gamma (X)=\{x\in E-X:X\cup x\in \mathcal{%
F}\}$ be the set of \textit{feasible continuations }of $X$. It is easy to
see that an accessible set system $(E,\mathcal{F})$ satisfies the interval
property without upper bounds if and only if for\ any\ $X,Y\in \mathcal{F},\
X\subseteq Y$\ implies\ $\Gamma (X)\cap (E-Y)\subseteq \Gamma (Y)$.

The maximal feasible subset of set $X\subseteq E$ is called a \textit{basis
of }$X$. \smallskip Clearly, by $(ii)$, there is only one basis for each
set. It will be\textit{\ }denoted by \-$\mathcal{B}(X)$.

\begin{definition}
For any $k\leq |E|$ the $k$\textit{-truncation }of a set system $(E,\mathcal{%
F})$ is a new set system defined by 
\[
\mathcal{F}_{k}=\{X\in \mathcal{F}:\left| X\right| \leq k\}.
\]
\end{definition}

If $(E,\mathcal{F)}$ is an antimatroid, then $(E,\mathcal{F}_{k})$ is the $k$%
\textit{-truncated antimatroid }\cite{BF}.

The \textit{rank} of a set $X\subseteq E$ is defined as $\varrho (X)=\max
\{\left| Y\right| :(Y\in \mathcal{F})\wedge (Y\subseteq X)\}$, the rank of
the set system $(E,\mathcal{F})$ is defined as $\varrho (\mathcal{F}%
)=\varrho (E).$ Notice, that every antimatroid $(E,\mathcal{F})$ is also a $%
k $-truncated antimatroid, where $k=\varrho (\mathcal{F})$.

A $k$-truncated antimatroid $(E,\mathcal{F})$ may not satisfy the interval
property without upper bounds, but it does satisfy the following condition: 
\begin{equation}
if\ X,Y\in \mathcal{F}_{k-1}\ and\ X\subseteq Y,\ then\ x\in E-Y,X\cup
x\in \mathcal{F}\ imply\ Y\cup x\in \mathcal{F}.  \label{trunc_isot}
\end{equation}

A set system $(E,\mathcal{F})$ has the $k$\textit{-truncated interval
property without upper bounds }if it satisfies (\ref{trunc_isot}).

\begin{theorem}
\cite{KL}An accessible set system $(E,\mathcal{F})$ of rank $k$ is a $k$%
-truncated antimatroid if and only if it satisfies the $k$-truncated
interval property without upper bounds.
\end{theorem}

In this paper we consider quasi-concave set functions on truncated
antimatroids.

\begin{definition}
A set function $F$ defined on a set system $(E,\mathcal{F})$ is
quasi-concave if for\ each\ $X,Y\in \mathcal{F}$,\ and\ for\ any\ maximal
feasible subset $Z$ of $X\cap Y$%
\begin{equation}
F(Z)\geq \min \{F(X),F(Y)\}.  \label{q_con}
\end{equation}
\end{definition}

Originally, these functions were considered \cite{Malishevski} on the
Boolean $2^{E}$, where the inequality (\ref{q_con}) turns into the following
condition 
\[
for\ each\ X,Y\subset E,\ F(X\cap Y)\geq \min \{F(X),F(Y)\}. 
\]
For this case, the correspondence between quasi-concave set functions and
monotone linkage functions were established in \cite{Mullat}.

A function $\pi :E\times 2^{E}\rightarrow \mathbf{R}$ is called \textit{a
monotone linkage function} if 
\begin{equation}
for\ all\ X,Y\subseteq E\ and\ x\in E,\ X\subseteq Y\ implies\ \pi (x,X)\geq
\pi (x,Y).  \label{three}
\end{equation}

Consider $F:2^{E}\rightarrow \mathbf{R}$ defined for each $X\subset E$%
\begin{equation}
F(X)=\min_{x\in E-X}\pi (x,X).  \label{four}
\end{equation}

It was shown \cite{Malishevski}, that $F$ is quasi-concave, and, moreover,
for every quasi-concave function $F$ there exists a monotone linkage
function $\pi $, which determines $F$ in accordance with (\ref{four}).

\smallskip In this work we extend these results to truncated antimatroids.
The family of feasible sets $\mathcal{F}$ of a truncated antimatroid $(E,%
\mathcal{F})$ forms a meet semilattice $L_{\mathcal{F}}$, with the lattice
operation: 
\[
X\wedge Y=\mathcal{B}(X\cap Y).
\]
Hence, for truncated antimatroids the inequality (\ref{q_con}) is converted
to the inequality 
\[
\ F(X\wedge Y)\geq \min \{F(X),F(Y)\}
\]
for\ each\ $X,Y\in L_{\mathcal{F}}$.

\section{Main results}

The following theorem characterizes quasi-concave functions defined on $k$%
-truncated antimatroids. Note, that in fact, we consider the functions
defined only on $\mathcal{F}_{k-1}$.

\begin{theorem}
\label{T-2}A set function $F$ defined on a $k$-truncated antimatroid $(E,%
\mathcal{F})$ is quasi-concave if and only if there exist a monotone linkage
function $\pi $ such that for each $X\in \mathcal{F}_{k-1}$%
\begin{equation}
F(X)=\min_{x\in \Gamma (X)}\pi (x,X).  \label{F_by_pi}
\end{equation}
\end{theorem}

%TCIMACRO{\TeXButton{Proof}{\proof} }
%BeginExpansion
%\proof%
%EndExpansion
\setlength {\parindent}{0.0cm}\textbf{Proof.}
Let a set function $F$ defined as a minimum of a monotone linkage function $%
\pi $. Note, that since for any antimatroid the operator $\Gamma $ is
not-empty for each $X\in \mathcal{F}_{k-1}$, the definition (\ref{F_by_pi})
is correct. To prove that the function $F$ is quasi-concave on $\mathcal{F}%
_{k-1}$, first note that
\setlength
{\parindent}{3.45ex} 
\begin{equation}
for\ each\ X\subset E,\ \Gamma (\mathcal{B}(X))\subseteq E-X,  \label{two}
\end{equation}
which immediately follows from the definition of basis.

Since $F(X\wedge Y)=\min_{x\in \Gamma (X\wedge Y)}\pi
(x,X\wedge Y)$ there is $x^{*}\in \Gamma (X\wedge Y)$ such that $F(X\wedge
Y)=\pi (x^{*},X\wedge Y)$. Then, by (\ref{two}), $x^{*}\in E-(X\cap Y)$,
i.e., either $x^{*}\in E-X$ or $x^{*}\in E-Y$. Without loss of generality,
assume that $x^{*}\in E-X$. Thus $X\wedge Y\subseteq X$, and $x^{*}\in E-X$,
and $x^{*}\in \Gamma (X\wedge Y)$, that accordingly to (\ref{trunc_isot})
implies $x^{*}\in \Gamma (X)$. Finally, 
\[
F(X\wedge Y)=\pi (x^{*},X\wedge Y)\geq \pi (x^{*},X)\geq \min_{x\in
\Gamma (X)}\pi (x,X)=F(X)\geq \min \{F(X),F(Y)\}. 
\]

To extend this function to the whole truncated antimatroid $(E,\mathcal{F})$
we can define $F(X)$ for each maximal $X$, i.e., for $|X|=k$, as $F(X)=%
\min_{(x,X)}\pi (x,X)$. It is easy to check that this extension
is quasi-concave too.

Conversely, let we have a quasi-concave set function $F$. Define the
function 
\begin{equation}
\pi _{F}(x,X)=\left\{ 
\begin{array}{ll}
\max_{A\in [X,E-x]_{\mathcal{F}_{k-1}}}F(A), & x\notin X\ and\
[X,E-x]_{\mathcal{F}_{k-1}}\neq \emptyset \\ 
\min_{A\in \mathcal{F}_{k-1}}F(A), & otherwise
\end{array}
\right. .  \label{pi_F}
\end{equation}

The function $\pi _{F}$ is monotone. Indeed, if $x\in E-Y$ and $[Y,E-x]_{%
\mathcal{F}_{k-1}}\neq \emptyset $, then $X\subseteq Y$ implies 
\[
\ \pi (x,X)=\max_{A\in [X,E-x]_{\mathcal{F}_{k-1}}}F(A)\geq 
\max_{A\in [Y,E-x]_{\mathcal{F}_{k-1}}}F(A)=\pi (x,Y). 
\]
It is easy to verify the remaining cases.

Let us denote $G(X)=\min_{x\in \Gamma (X)}\pi _{F}(x,X)$, and
prove that $F=G$ on $\mathcal{F}_{k-1}$.

Now 
\[
G(X)=\min_{x\in \Gamma (X)}\pi _{F}(x,X)=\pi _{F}(x^{*},X)=%
\max_{A\in [X,E-x^{*}]_{\mathcal{F}_{k-1}}}F(A)\geq F(X). 
\]

On the other hand, 

\[
G(X)=\min_{x\in \Gamma (X)}\pi _{F}(x,X)=%
\min_{x\in \Gamma (X)}F(A^{x}),
\]
where $A^{x}$ is a set from $%
[X,E-x]_{\mathcal{F}_{k-1}}$ on which the value of the function $F$ is
maximal, i.e., 

\[
A^{x}=\arg\max_{A\in [X,E-x]_{\mathcal{F}_{k-1}}}F(A). 
\]
From quasi-concavity of $F$ follows that 
\[
\min_{x\in \Gamma (X)}F(A^{x})\leq F(\wedge_{x\in \Gamma (X)%
}A^{x}).
\]
So, $G(X)\leq F(\wedge_{x\in \Gamma (X)}A^{x})$.

It remains to prove, that for all $X\in \mathcal{F}_{k-1}$, $X=\wedge_{%
x\in \Gamma (X)}A^{x}$, where the set $A^{x}\in [X,E-x]_{\mathcal{F}_{k-1}}$%
.

Denote, $Y=\wedge_{x\in \Gamma (X)}A^{x}$. For each $x\in
\Gamma (X)$, $X\subseteq A^{x}$, and consequently $X\subseteq Y$. Assume,
that $X\subset Y$, then by definition $(A2)$ there exists an element $y\in
Y-X$ such that $X\cup y\in \mathcal{F}$, i.e., $y\in Y\cap \Gamma (X)$. On
the other hand, 
\[
Y=\wedge_{x\in \Gamma (X)}A^{x}\subseteq \cap_{x\in
\Gamma (X)}A^{x}\subseteq E-\Gamma (X). 
\]
This contradiction proves that $X=Y$.

Therefore, $G(X)\leq F(X)$, and, hence, $F=G$, i.e. $F(X)=\min_{x\in
\Gamma (X)}\pi _{F}(x,X)$, where $\pi _{F}$ is a monotone linkage
function. %
%TCIMACRO{\TeXButton{End Proof}{\endproof}}
%BeginExpansion
%\endproof%
%EndExpansion
\rule{2mm}{2mm} 

\smallskip Thus, we proved that each quasi-concave function $F$ determines a monotone
linkage function $\pi _{F}$, and the set function defined as the minimum of
this monotone linkage function $\pi _{F}$ coincides with the original
function $F$. A weaker property holds for the linkage functions.

\begin{theorem}
\label{null_pi}Let $F(X)=\min_{x\in \Gamma (X)}\pi _{F}(x,X)$
for a monotone linkage function $\pi $ on a $k$-truncated antimatroid $(E,%
\mathcal{F})$. Then $\pi _{F}|_{\mathcal{F}_{k-1}}\leq \pi |_{\mathcal{F}%
_{k-1}}$, i.e., for any $X\in \mathcal{F}_{k-1}$ and $x\in \Gamma (X)\ $%
\[
\pi _{F}(x,X)\leq \pi (x,X),
\]
where $\pi _{F}$ is defined by (\ref{pi_F}).
\end{theorem}

%TCIMACRO{\TeXButton{Proof}{\proof} }
%BeginExpansion
%\proof%
%EndExpansion
\setlength {\parindent}{0.0cm}\textbf{Proof.}
For any $X\in \mathcal{F}_{k-1}$ and $x\in \Gamma (X)$%
\setlength
{\parindent}{3.45ex}
\[
\pi _{F}(x,X)=\max_{A\in [X,E-x]_{\mathcal{F}_{k-1}}}%
F(A)=F(A^{*})=\min_{a\in \Gamma (A^{*})}\pi (a,A^{*})\leq \pi
(x,A^{*}). 
\]

The last inequality follows from the $k$-truncated interval property without
upper bounds. Indeed, $X\subseteq A^{*}$ and $x\notin A^{*}$, then $x\in
\Gamma (X)$ implies $x\in \Gamma (A^{*})$.

Now, from monotonicity of the function $\pi $ we have $\pi (x,A^{*})\leq \pi
(x,X)$, that finishes the proof. %
%TCIMACRO{\TeXButton{End Proof}{\endproof}}
%BeginExpansion
%\endproof%
%EndExpansion
\rule{2mm}{2mm}

\smallskip Consider the following example to see that these two functions $\pi $ and $%
\pi _{F}$ may be not equal. For example, let $E=\{1,2\}$, $\mathcal{F}=2^{E}$%
, and 
\[
\pi (x,X)=\left\{ 
\begin{array}{ll}
2, & x=2\mathit{\ }and\ X=\emptyset \\ 
1, & otherwise.
\end{array}
\right. 
\]
Then the function $F(X)=\min_{x\in \Gamma (X)}\pi (x,X)$ is
equal to $1$ for all $X\subset E$, and $\pi _{F}$ equals for $1$ for each
pair $(x,X)\in E\times 2^{E}$, i.e., $\pi _{F}$ $\neq \pi $.

Now let us define more exactly the structure of the set of monotone linkage
functions.

\begin{theorem}
Let $(E,\mathcal{F})$ be a set system of rank $k$, where the set of feasible
continuations\textit{\ }of $X$ is not empty for each $X\in \mathcal{F}_{k-1}$%
, and let $\pi _{1\ }$and $\pi _{2}$ define (by (\ref{F_by_pi})) the
same set function $F$ on $\mathcal{F}_{k-1}$. Then the function 
\[
\pi =\min \{\pi _{1\ },\pi _{2}\}
\]

is a monotone linkage function, and it determines the same function $F$ on $%
\mathcal{F}_{k-1}$.
\end{theorem}

%TCIMACRO{\TeXButton{Proof}{\proof}}
%BeginExpansion
%\proof%
%EndExpansion
\setlength {\parindent}{0.0cm}\textbf{Proof.}
At first, prove that $\pi $ is a monotone linkage function. Indeed, consider
a pair $X\subseteq Y$. Suppose, without loss of generality, that $\min \{\pi
_{1}(x,X),\pi _{2}(x,X)\}=\pi _{1}(x,X)$. Now, 
\setlength
{\parindent}{3.45ex}
\[
\pi (x,X)=\min \{\pi _{1}(x,X),\pi _{2}(x,X)\}=\pi _{1}(x,X)\geq 
\]
\[
\geq \pi _{1}(x,Y)\geq \min \{\pi _{1}(x,Y),\pi _{2}(x,Y)\}=\pi (x,Y)
\]

To complete the proof, we show that 
\[
\min_{x\in \Gamma (X)}\pi (x,X)=\pi (x^{*},X)=\min \{\pi
_{1}(x^{*},X),\pi _{2}(x^{*},X)\}\geq 
\]
\[
\geq \min (\min_{x\in \Gamma (X)}\pi _{1}(x,X),\min_{x\in
\Gamma (X)}\pi _{2}(x,X))=F(X),
\]

and on the other hand, 
\[
F(X)=\min_{x\in \Gamma (X)}\pi _{1}(x,X)=\pi _{1}(x^{\#},X)\geq
\pi (x^{\#},X)\geq \min_{x\in \Gamma (X)}\pi (x,X). 
\]
%TCIMACRO{\TeXButton{End Proof}{\endproof}}
%BeginExpansion
%\endproof%
%EndExpansion
\rule{2mm}{2mm}

\smallskip Thus, the set of monotone linkage functions, defining a set function $F$ on
a truncated antimatroid, forms a semilattice with the following lattice
operation: 
\[
\pi _{1}\wedge \pi _{2}=\min \{\pi _{1\ },\pi _{2}\}, 
\]
where by Theorem \ref{null_pi} the function $\pi _{F}$ is a null of this
semilattice.

The following theorem demonstrates the necessity of interval property for
the above established correspondence between quasi-concave set functions and
monotone linkage functions.

\begin{theorem}
\smallskip Let $(E,\mathcal{F})$ be an accessible set system of rank $k$. If
the set of feasible continuations\textit{\ }of $X$ is not empty for each $%
X\in \mathcal{F}_{k-1}$, then the following statements are equivalent

$(i)$ $(E,\mathcal{F})$ is a $k$-truncated antimatroid

$(ii)$ the function $F=\min_{x\in \Gamma (X)}\pi (x,X)$ is
quasi-concave for every monotone linkage function $\pi $.
\end{theorem}

%TCIMACRO{\TeXButton{Proof}{\proof}}
%BeginExpansion
%\proof%
%EndExpansion
\setlength {\parindent}{0.0cm}\textbf{Proof.}
Since the one direction is proved (see Theorem \ref{T-2}), assume that the
set system $(E,\mathcal{F})$ is not $k$-truncated antimatroid, i.e., there
exist $A,B\in \mathcal{F}_{k-1}$ such that $A\subset B$, and there is $a\in
E-B$ such that $A\cup a\in \mathcal{F}$ and $B\cup a\notin \mathcal{F}$.
Define the linkage function 
\setlength
{\parindent}{3.45ex}
\[
\pi (x,X)=\left\{ 
\begin{array}{ll}
0, & x\in X \\ 
1, & x=a\mathit{\ }and\ A\subseteq X\subseteq E-a \\ 
2, & otherwise
\end{array}
\right. 
\]

It is easy to check that $\pi $ is monotone.

Here, $F(A)=1$, $F(A\cup a)=$ $F(B)=2$. Since $(A\cup a)\cap B=A$, we have 
\[
F((A\cup a)\cap B)<\min \{F(A\cup a),F(B)\}, 
\]
i.e., $F$ is not quasi-concave. 
%TCIMACRO{\TeXButton{End Proof}{\endproof}}
%BeginExpansion
%\endproof%
%EndExpansion
\rule{2mm}{2mm}

\section{\protect\smallskip Conclusions}

In this article, we discuss the duality between quasi-concave set functions
and monotone linkage functions. It is shown that each quasi-concave function 
$F$, defined on an antimatroid, determines a semilattice of monotone linkage
functions each of them defines the set function $F$, and the null of this
semilattice is the function $\pi _{F}$ constructed from the function $F$.

As the directions for future research we see the extension of the duality to
other families of sets such as convex geometries, interval greedoids and
more general set families.

\end{document}